# Real-Time Topology Detection and State Estimation in Distribution Systems Using Micro-PMU And Smart Meter Data

Zahra Soltani, *Student Member, IEEE*, and Mojdeh Khorsand, *Member, IEEE*

***Abstract*— Distribution network topology detection and state estimation in real-time are critical for modern distribution systems management and control. However, number of sensors in distribution networks are limited and communication links between switch devices and distribution management system are not well-established. In this regard, this paper proposes mixed-integer quadratic programming (MIQP) formulations to determine the topology of distribution network and estimate distribution system states simultaneously using micro-phasor measurement units (micro-PMUs) and smart meter data. Two approaches based on AC optimal power flow are proposed and analyzed: (i) polar power-voltage (PPV) formulation, and (ii) rectangular current-voltage (RIV) formulation. The proposed models include convex objective function while constraints are linearized using first-order approximation of Taylor series and Big M method. The proposed models can identify multiple simultaneous switching actions at each time instant and different topology configurations including radial and meshed networks. Only measurement data at each time interval is needed to identify topology and states of the system correctly. The proposed models are tested on a modified IEEE-33 bus system with realistic load data from Pecan Street Inc. database. The results confirm that both models can identify system topology and states with remarkable accuracy in real-time, while RIV model outperforms PPV model.**

***Index Terms*- Topology detection, state estimation, distribution system, micro-PMU, smart meter, mixed-integer quadratic programming, rectangular current-voltage formulation, polar power-voltage formulation, AC optimal power flow.**

## Nomenclature

*Sets and Indices*

| | |
|---|---|
| $a, c$ | Index for bus, $a, c \in \psi$ |
| $r$ | Index for micro-PMU, $r \in \Omega$ |
| $y$ | Index for measurement |
| $e$ | Index for switch, $e \in \Upsilon$ |

*Parameters and Constants*

| | |
|---|---|
| $g_{ac}$ | Conductance of line connecting bus $a$ and bus $c$ |
| $b_{ac}$ | Susceptance of line connecting bus $a$ and bus $c$ |
| $E_r^M$ | Measured voltage magnitude by micro-PMU $r$ |
| $E_r^P$ | Actual voltage magnitude at micro-PMU $r$ location |
| $\varphi_r^M$ | Measured voltage angle by micro-PMU $r$ |
| $\varphi_r^P$ | Actual voltage angle at micro-PMU $r$ location |
| $w_y$ | Weight of measurement $y$ |
| $M_i$ | A large positive number; $i = 1, 2, 3, or\ 4$ |

*Variables*

| | |
|---|---|
| $P_{a,c}^L$ | Active power flow from bus $a$ to bus $c$ |
| $Q_{a,c}^L$ | Reactive power flow from bus $a$ to bus $c$ |
| $E_a$ | Voltage magnitude at bus $a$ |
| $\varphi_a$ | Voltage angle at bus $a$ |
| $P_{k,a}^D$ | Active demand of load $k$ at bus $a$ |
| $Q_{k,a}^D$ | Reactive demand of load $k$ at bus $a$ |
| $\vartheta_e^l$ | Status of switch $e$ |
| $P_{i,a}^G$ | Active power of generator $i$ at bus $a$ |
| $Q_{i,a}^G$ | Reactive power of generator $i$ at bus $a$ |
| $I_{a,c}^r$ | Real part of current from bus $a$ to $c$ |
| $I_{a,c}^{im}$ | Imaginary part of current from bus $a$ to $c$ |
| $E_{a,c}^r$ | Real part of voltage at bus $a$ |
| $E_{a,c}^{im}$ | Imaginary part of voltage at bus $a$ |
| $I_a^r$ | Real part of current injection at bus $a$ |
| $I_a^{im}$ | Imaginary part of current injection at bus $a$ |
| $\lambda_{a,c}^L$ | Auxiliary variable |
| $\chi_{a,c}^L$ | Auxiliary variable |

## I. Introduction

THE information of distribution network topology and system states are crucial for real-time (RT) operation and control of distribution systems, e.g., volt-Var control, especially for systems with high penetration of distributed energy resources (DERs). However, communication links are not installed for majority of switch devices in the distribution networks, which makes it difficult to maintain an updated network topology information in distribution management systems. Moreover, only limited number of sensors are installed in the distribution networks, which provide incomplete observability of the system for the distribution system operator. Also, integration of DERs may result in more reconfiguration and switching actions in the distribution system. Thus, an efficient distribution system topology processor and state estimation tool is critical for success of distribution management systems.

For enhanced reliability, modern distribution systems for urban areas are often designed with a loosely meshed or looped connection between feeders or substations. Even though the system may be operated radially, the loop provides more than one point of interconnection, improves efficiency and

The research is funded by the Department of Energy (DOE) Advanced Research Projects Agency – Energy (ARPA-E) under OPEN 2018 program.

Zahra Soltani and Mojdeh Khorsand are with the School of Electrical, Computer, and Energy Engineering, Arizona State University, Tempe, AZ 85281 USA (e-mail: zsoltani@asu.edu; mojdeh.khorsand@asu.edu).

reliability, and prevents transmission fault currents from flowing across the distribution system and damaging equipment while reducing load shedding. Moreover, meshed distribution systems exist in many metropolitan areas. Also, networked microgrids are emerging within distribution systems. Recent research has evidenced that weakly-meshed operations may yield significant benefits including improvements in balancing power, losses, voltage profiles, and higher hosting capacity for distributed generation (DG) [1]-[2]. The transformation of distribution systems from passive to active networks with DERs and meshed or weakly-meshed structures highlights the need for an efficient topology processor. In [3], a model is provided to optimize the sensor placement for topology identification. For a particular location of sensors, this model gives the confidence level of identifying changes in switches status. Residual error obtained from the state estimation is used to identify network topology [4]. A recursive Bayesian approach is employed in [5] to perform state estimation for all possible topologies and identify the topology with highest probability as the correct network topology. However, the algorithm presented in [5] is not computationally efficient. The reason is that, for any possible topology configurations of a distribution network, this method performs state estimation, and then chooses the topology with the highest probability. A topology error detection method based on state estimation is proposed in [6], where the circuit breakers statuses are considered as state variables and telemetered statuses of circuit breakers are incorporated into the model. However, the method proposed in [6] may not be applicable to the distribution networks with limited number of telemetered switches. Moreover, data-driven approaches for topology processor have been proposed in [7]–[14]. Voltage correlation analysis is utilized in [7] to detect the distribution network topology using graph theory. A graph learning approach is proposed in [8] to reconstruct feeder topologies in distribution systems based on nodal voltage measurements. Smart meters and micro-phasor measurement units (micro-PMUs) have gained reputation in monitoring of power distribution systems [10]. Micro-PMUs provide synchronized measurements of voltage and current phasors [11]. Using the smart meter data for building voltage covariance, a maximum a-posteriori probability method is proposed in [9] to identify topology of a distribution system. Time-series signature verification method for identifying the topology of distribution network based on measured voltages by micro-PMUs has been initially proposed in [12]-[13], which assumes the same resistance to reactance ratio for all electric grid lines. This method is further developed in [14], in which based on the prior information of switch statuses, a library of signatures is calculated to obtain possible topology configurations. Then, the change in the voltage time series measured by micro-PMUs are compared with the obtained library to detect change in the topology of distribution system. The main drawback of [12]–[14] is that the authors assume that the topology change may occur due to only one switching action at each time. Also, the prior information of switch statuses and prior voltage measured by micro-MPUs are needed to identify the network topology. In this regard, if the load variation is increased, or the prior status of switches is obtained wrongly, the topology may not be identified correctly. Furthermore, this method is dependent to three parameter tunings. In [15], a single-shot mixed-integer quadratic programming (MIQP) problem is proposed based on DC power flow assumptions to obtain the circuit breaker statuses at substations. However, the assumptions of DC power flow model are not appropriate for the topology processor in the distribution networks. The distribution network topology processing and state estimation problem is a mixed-integer nonlinear programming (MINLP) problem due to binary variables associated with status of switches and nonlinear AC power flow equations.

In this paper, two MIQP formulations are proposed to determine the topology of distribution network and estimate system states simultaneously using micro-PMUs and smart meters data. The proposed models are able to identify different topology configurations including radial, looped, and meshed networks. The proposed MIQP approaches are based on two AC optimal power flow models, (i) polar power-voltage (PPV) formulation and (ii) rectangular current-voltage (RIV) formulation, which are linearized using iterative first-order approximation of Taylor series. The performance of these approaches is compared under load's variability and measurement noises. In order to eliminate nonlinearity due to inclusion of binary variables associated with status of switches, the big M technique, which has been used in the authors' prior work for transmission switching is leveraged [16]–[20]. The proposed AC optimal power flow models include linear constraints and convex objective functions, which can obtain the global optimal solution via optimization solvers utilizing the branch and bound algorithm to solve MIQP problems. The proposed approaches are able to identify multiple simultaneous switching actions at each time instant without information of switch statuses in prior time intervals. The proposed models are single-shot optimization problems, i.e., they only require measurement data at each time snapshot to identify the topology of the system and estimate system states accurately.

The rest of the paper is organized as follows. Sections II and III show PPV-based and RIV-based topology detection and state estimation formulation in the distribution network, respectively. In Section IV, case studies and simulation results are provided. Section V presents conclusion.

## II. PPV-BASED TOPOLOGY DETECTION AND STATE ESTIMATION MODEL IN DISTRIBUTION NETWORK

In this section, PPV-based simultaneous topology detection and state estimation model in a distribution system using micro-PMUs and smart meters data is discussed. First, the nonlinear PPV-based model is explained. Second, the formulation of proposed MIQP-PPV-model is presented.

### *A. PPV-based Topology Detection and State Estimation Formulation*

The nonlinear AC power flow equations can be formulated in various forms including PPV model or RIV model. In this section, the PPV-based topology detection and state estimation problem in distribution networks is formulated, which is valid for meshed, looped, or radial topology structures. Assume a distribution network with set of buses $\psi = \{1,2,\dots,O\}$ and set of lines $\Phi = \{1,2,\dots,P\}$. Set of micro-PMUs is represented by $\Omega = \{1,2,\dots,R\}$. For the line $p \in \Phi$, which connects bus $a \in \psi$ to bus $c \in \psi$ and is always energized, i.e., it is non-switchable, the nonlinear active and reactive AC power flow equations are defined using (1.a)-(1.b) [21]-[22].

$$P^L_{a,c} = {E_a}^2 g_{ac} - E_a E_c g_{ac}\cos(\varphi_a - \varphi_c) - E_a E_c b_{ac}\sin(\varphi_a - \varphi_c),\ \forall (a,c) \in \Phi \tag{1.a}$$

$$Q^L_{a,c} = -{E_a}^2 b_{ac} + E_a E_c\, b_{ac}\cos(\varphi_a - \varphi_c) - E_a E_c\, g_{ac}\sin(\varphi_a - \varphi_c),\ \forall (a,c) \in \Phi \tag{1.b}$$

Let $\Upsilon = \{1,2,\dots,E\}$ be the set of switches in a distribution system. The active and reactive power flow in the line $p \in \Phi$ equipped with a switch device $e \in \Upsilon$ is modeled by including binary variable $\vartheta^l_e$ in (1.a)-(1.b) as follows:

$$P^L_{a,c} = \vartheta^l_e\big({E_a}^2 g_{ac} - E_a E_c g_{ac}\cos(\varphi_a - \varphi_c) - E_a E_c b_{ac}\sin(\varphi_a - \varphi_c)\big),\forall (a,c) \in \Phi \tag{1.c}$$

$$Q^L_{a,c} = \vartheta^l_e\big(-{E_a}^2 b_{ac} + E_a E_c\, b_{ac}\cos(\varphi_a - \varphi_c) - E_a E_c\, g_{ac}\sin(\varphi_a - \varphi_c)\big),\forall (a,c) \in \Phi \tag{1.d}$$

where $\vartheta^l_e = 0$ indicates the switch is open and the line is disconnected, while $\vartheta^l_e = 1$ implies the switch is closed and the line is energized.

The active and reactive power balance constraints at bus $a \in \psi$ in a distribution network are given by:

$$\sum_{\forall i \epsilon i(a)} P^G_{i,a} = \sum_{\substack{\forall c \in \psi(a) \\ c \neq a}} P^L_{a,c} + \sum_{\forall k \epsilon K(a)} P^D_{k,a},\ \forall a \in \psi \tag{1.e}$$

$$\sum_{\forall i \epsilon i(a)} Q^G_{i,a} = \sum_{\substack{\forall c \in \psi(a) \\ c \neq a}} Q^L_{a,c} + \sum_{\forall k \epsilon K(a)} Q^D_{k,a},\ \forall a \in \psi \tag{1.f}$$

The synchronized voltage magnitude and phase angle measurements provided by the micro-PMUs not only improve the real-time monitoring of distribution system, but also provide direct measurement of system states [23]. However, the number of micro-PMUs is limited to only few in distribution systems. To evaluate micro-PMU measurement noise, the total vector error (TVE) index is used [14]. TVE is expressed as normalized value of the difference between actual and measured phasor values. The micro-PMU voltage phasor measurement $r \in \Omega$ can be modeled by (1-g)-(1.h).

$$E^M_r = E^P_r + \mu_r \tag{1.g}$$

$$\varphi^M_r = \varphi^P_r + \gamma_r \tag{1.h}$$

where $\mu_r$ and $\gamma_r$ are Gaussian noises with respect to TVE index. The PPV-based topology detection and state estimation formulation in the distribution system is proposed as follows:

$$Min \sum_{y=1}^{Y} w_y (h_y(f) - z_y)^2 \tag{1.i}$$

$$Subject\ to\ (1.a) - (1.f).$$

where $z_y$ is measurement value $y$, $f$ is a vector of the system states including $E$ and $\varphi$, and $h_y(f)$ is nonlinear function of system states related to the measurements in a distribution network, which include substation, smart meter, and micro-PMU measurements. The vector $\Lambda = \{\vartheta^l_1, \vartheta^l_2, \dots, \vartheta^l_E\}$ represents the network topology.

### *B. Proposed MIQP-PPV-based Topology Detection and State Estimation Formulation*

The PPV-based distribution network topology detection and state estimation problem in (1) is a MINLP problem. The nonlinear terms are product of binary variable $\vartheta^l_{e,x}$ and continuous variables as well as the nonlinear active and reactive AC power flow equations. Such problem can be solved using nonlinear algorithms, which may diverge or obtain local optimal solutions. A MIQP model based on DC power flow is proposed in [15] to determine the breaker statuses at substations. However, DC power flow model is not suitable for the topology processor in the distribution networks. To cope with such challenges, a MIQP formulation based on a linearized PPV (MIQP-PPV-based) AC power flow model is proposed in this paper to determine the topology and states of a distribution system using micro-PMUs and smart meters measurements. To this end, first, the linear approximations of nonlinear active and reactive AC power flow constraints in (1) are proposed using the iterative first-order approximation of Taylor series, which are defined in (2.a)-(2.b).

$$P^L_{a,c} = g_{ac}\vartheta^l_e\Big[-{\underline{E_{a,it-1}}}^2 + 2\underline{E_{a,it-1}}E_a - \underline{E_{c,it-1}}\cos\left(\underline{\varphi_{ac,it-1}}\right) E_a + \underline{E_{a,it-1}E_{c,it-1}}\cos\left(\underline{\varphi_{ac,it-1}}\right) - \underline{E_{a,it-1}}\cos\left(\underline{\varphi_{ac,it-1}}\right) E_c + \varphi_{ac}\underline{E_{a,it-1}E_{c,it-1}}\sin\left(\underline{\varphi_{ac,it-1}}\right) - \underline{\varphi_{ac,it-1}E_{a,it-1}E_{c,it-1}}\sin\left(\underline{\varphi_{ac,it-1}}\right)\Big] - b_{ac}\vartheta^l_e\Big[\sin\left(\underline{\varphi_{ac,it-1}}\right)\underline{E_{c,it-1}}E_a - \underline{E_{a,it-1}}\,\underline{E_{c,it-1}}\sin\left(\underline{\varphi_{ac,it-1}}\right) + \underline{E_{a,it-1}}\sin\left(\underline{\varphi_{ac,it-1}}\right)E_c + \varphi_{ac}\underline{E_{a,it-1}}\,\underline{E_{c,it-1}}\cos\left(\underline{\varphi_{ac,it-1}}\right) - \underline{\varphi_{ac,it-1}E_{a,it-1}}\,\underline{E_{c,it-1}}\cos\left(\underline{\varphi_{ac,it-1}}\right)\Big],\ \forall (a,c) \in \Phi \tag{2.a}$$

$$Q^L_{a,c} = b_{ac}\vartheta^l_e\Big[{\underline{E_{a,it-1}}}^2 - 2\underline{E_{a,it-1}}E_a + \underline{E_{c,it-1}}\cos\left(\underline{\varphi_{ac,it-1}}\right) E_a - \underline{E_{a,it-1}E_{c,it-1}}\cos\left(\underline{\varphi_{ac,it-1}}\right) + \underline{E_{a,it-1}}\cos\left(\underline{\varphi_{ac,it-1}}\right) E_{c,x} - \varphi_{ac}\underline{E_{a,it-1}E_{c,it-1}}\sin\left(\underline{\varphi_{ac,it-1}}\right) + \underline{\varphi_{ac,it-1}E_{a,it-1}E_{c,it-1}}\sin\left(\underline{\varphi_{ac,it-1}}\right)\Big] - g_{ac}\vartheta^l_e\Big[\sin\left(\underline{\varphi_{ac,it-1}}\right)\underline{E_{c,it-1}}E_a - \underline{E_{a,it-1}}\,\underline{E_{c,it-1}}\sin\left(\underline{\varphi_{ac,it-1}}\right) + \underline{E_{a,it-1}}\sin\left(\underline{\varphi_{ac,it-1}}\right)E_c + \varphi_{ac}\underline{E_{a,it-1}E_{c,it-1}}\cos\left(\underline{\varphi_{ac,it-1}}\right) - \underline{\varphi_{ac,it-1}E_{a,it-1}}\,\underline{E_{c,it-1}}\cos\left(\underline{\varphi_{ac,it-1}}\right)\Big],\ \forall (a,c) \in \Phi \tag{2.b}$$

where $\underline{E_{a,it-1}}$, $\underline{E_{c,it-1}}$, and $\underline{\varphi_{ac,it-1}}$ are updated in each iteration of the proposed MIQP-PPV-based model based on solution of previous iteration. However, (2.a)-(2.b) are still nonlinear as a result of multiplication of $\vartheta^l_e$ and continuous variables $E_a$, $E_c$, $\varphi_a$, and $\varphi_c$. To eliminate such nonlinearity, big $M$ technique, which has been used in authors' prior work for topology control, is leveraged [16]–[20]. The nonlinear equations (2.a)-(2.b) are linearized using big $M$ method as follows:

$$-M_1(1 - \vartheta^l_e) \le \chi^L_{a,c} - P^L_{a,c} \le M_1(1 - \vartheta^l_e) \tag{2.c}$$

$$-M_1\vartheta^l_e \le P^L_{a,c} \le M_1\vartheta^l_e \tag{2.d}$$

$$-M_2(1 - \vartheta^l_e) \le \lambda^L_{a,c} - Q^L_{a,c} \le M_2(1 - \vartheta^l_e) \tag{2.e}$$

$$-M_2\vartheta_e^l \le Q_{a,c}^L \le M_2\vartheta_e^l \quad (2.f)$$

where $\chi_{a,c}^L$ and $\lambda_{a,c}^L$ are the linear parts of (2.a) and (2.b), i.e., right hand side of (2.a) and (2.b) without $\vartheta_e^l$ multiplier. Based on the proposed linear constraints, if $\vartheta_e^l$ is equal to zero, $P_{a,c,x}^L$ and $Q_{a,c,x}^L$ will become zero, and constraints (2.c) and (2.e) will not be binding.

The proposed joint PPV-based topology detection and state estimation formulation in the radial, looped, and meshed distribution system is formulated as follows:

$$Min \sum_{y=1}^{Y} w_y (h_y(f) - z_y)^2 \quad (2.g)$$

$Subject\ to$ (1.e)-(1.f), (2.a)-(2.f)

where $h_y(f)$ is the linear function of system states associated with the measurements. For non-switchable lines, multiplier $\vartheta_e^l$ is eliminated from (2.a)-(2.b), i.e., $\vartheta_e^l = 1$. For switchable lines, (2.c)-(2.f) are considered. The proposed model of (2.g) is MIQP with convex objective function and linear constraints.

## III. Proposed MIQP-RIV-based Topology Detection and State Estimation in Distribution Network

In this section, the linearized RIV-based distribution network topology detection and state estimation problem is proposed based on MIQP formulation. The current flow of non-switchable line $p \in \Phi$, which connects bus $a \in \psi$ to bus $c \in \psi$ can be obtained using linear constraints (3.a) and (3.b) [24].

$$I_{a,c}^r = g_{ac}(E_a^r - E_c^r) - b_{ac}(E_a^{im} - E_c^{im}), \forall (a,c) \in \Phi \quad (3.a)$$

$$I_{a,c}^{im} = b_{ac}(E_a^r - E_c^r) + g_{ac}(E_a^{im} - E_c^{im}), \forall (a,c) \in \Phi \quad (3.b)$$

If there is a switch device $e \in \Upsilon$ on the line $p \in \Phi$, equations (3.a)-(3.b) are modified by considering binary variable $\vartheta_e^l$ associated with status of the switch as follows:

$$I_{a,c}^r = \vartheta_e^l (g_{ac}(E_a^r - E_c^r) - b_{ac}(E_a^{im} - E_c^{im})), \forall (a,c) \in \Phi \quad (3.c)$$

$$I_{a,c}^{im} = \vartheta_e^l (b_{ac}(E_a^r - E_c^r) + g_{ac}(E_a^{im} - E_c^{im})), \forall (a,c) \in \Phi \quad (3.d)$$

Constraints (3.c)-(3.d) are nonlinear due to product of binary variable $\vartheta_e^l$ with continuous variables. In order to eliminate such nonlinearity, the big M method is utilized in this paper to linearize constraints (3.c)-(3.d) as follows:

$$-M_3(1 - \vartheta_e^l) \le I_{a,c}^r - \begin{bmatrix} g_{ac}(E_a^r - E_c^r) - \\ b_{ac}(E_a^{im} - E_c^{im}) \end{bmatrix} \quad (3.e)$$

$$I_{a,c}^r - \begin{bmatrix} g_{ac}(E_a^r - E_c^r) - \\ b_{ac}(E_a^{im} - E_c^{im}) \end{bmatrix} \le M_3(1 - \vartheta_e^l) \quad (3.f)$$

$$-M_3\vartheta_e^l \le I_{a,c}^r \le M_3\vartheta_e^l \quad (3.g)$$

$$-M_4(1 - \vartheta_e^l) \le I_{a,c}^{im} - \begin{bmatrix} b_{ac}(E_a^r - E_c^r) + \\ g_{ac}(E_a^{im} - E_c^{im}) \end{bmatrix} \quad (3.h)$$

$$I_{a,c}^{im} - \begin{bmatrix} b_{ac}(E_a^r - E_c^r) + \\ g_{ac}(E_a^{im} - E_c^{im}) \end{bmatrix} \le M_4(1 - \vartheta_e^l) \quad (3.i)$$

$$-M_4\vartheta_e^l \le I_{a,c}^{im} \le M_4\vartheta_e^l \quad (3.j)$$

The current injection constraints at bus $a \in \psi$ of a distribution system are formulated as (3.k)-(3.l).

$$I_a^r = \sum_{c \epsilon \delta(a)} I_{a,c}^r, \forall a \in \psi \quad (3.k)$$

$$I_a^{im} = \sum_{c \epsilon \delta(a)} I_{a,c}^{im}, \forall a \in \psi \quad (3.l)$$

The nonlinear active and reactive power injection constraints at bus $a \in \psi$ of the system are expressed in (3.m)-(3.n) [24].

$$\sum_{\forall i \epsilon i(a)} P_{i,a}^G - \sum_{\forall k \epsilon K(a)} P_{k,a}^D = E_a^r I_a^r + E_a^{im} I_a^{im}, \forall a \in \psi \quad (3.m)$$

$$\sum_{\forall i \epsilon i(a)} Q_{i,a}^G - \sum_{\forall k \epsilon K(a)} Q_{k,a}^D = E_{a,x}^{im} I_a^r - E_a^r I_a^{im}, \forall a \in \psi \quad (3.n)$$

The nonlinear active and reactive power injection constraints for bus $a \in \psi$ of the system are formulated as linear constraints (3.o) and (3.p) using iterative first-order approximation of Taylor series, respectively.

$$\sum_{\forall i \epsilon i(a)} P_{i,a}^G - \sum_{\forall k \epsilon K(a)} P_{k,a}^D = \underline{E_{a,it-1}^r} I_a^r + \underline{E_{a,it-1}^{im}} I_a^{im} + \underline{I_{a,t-1}^r} E_a^r + \underline{I_{a,it-1}^{im}} E_a^{im} - \underline{E_{a,it-1}^r \ I_{a,it-1}^r} - \underline{E_{a,it-1}^{im} \ I_{a,it-1}^{im}}, \forall a \in \psi \quad (3.o)$$

$$\sum_{\forall i \epsilon i(a)} Q_{i,a}^G - \sum_{\forall k \epsilon K(a)} Q_{k,a}^D = \underline{E_{a,it-1}^{im}} I_a^r - \underline{E_{a,it-1}^r} I_a^{im} + \underline{I_{a,it-1}^r} E_a^{im} - \underline{I_{a,it-1}^{im}} E_a^r - \underline{E_{a,it-1}^{im} \ I_{a,it-1}^r + E_{a,it-1}^r \ I_{a,it-1}^{im}}, \forall a \in \psi \quad (3.p)$$

where $\underline{E_{a,it-1}^r}$, $\underline{E_{a,it-1}^{im}}$, $\underline{I_{a,it-1}^r}$, and $\underline{I_{a,it-1}^{im}}$ are updated in each iteration of the proposed RIV-based model based on solution of the previous iteration. The proposed MIQP problem based on the linearized RIV (MIQP-RIV-based) model is formulated in (3.q) to identify topology and states of a distribution network.

$$Min \sum_{y=1}^{Y} w_y (h_y(f) - z_y)^2 \quad (3.q)$$

$Subject\ to$ (3.a)-(3.b), (3.e)-(3.l),(3.o)-(3.p)

where $f$ is a vector of the system states including $E_a^r$ and $E_a^{im}$ and $h_y(f)$ is the linear function of system states associated with the measurements. The proposed RIV-based model comprises convex objective function and linear constraints.

## IV. Proposed Iterative MIQP-based Topology Detection and State Estimation

Figure 1 shows flowchart of the simulation procedure for the proposed iterative MIQP-PPV-based and MIQP-RIV-based topology processor and state estimation models in distribution systems. The measurement data for two models include micro-PMUs, smart meters, and substation measurements. Flat start point (in PPV formulation, voltage magnitude = 1 and voltage angle = 0; in RIV formulation, real part of voltage=1 and imaginary part of voltage=0) is considered in the first iteration for all buses. Then, the proposed PPV-MIQP-based (given in (2.g)) and RIV-MIQP-based (given in (3.q)) models are solved to identify topology and states of the distribution system.

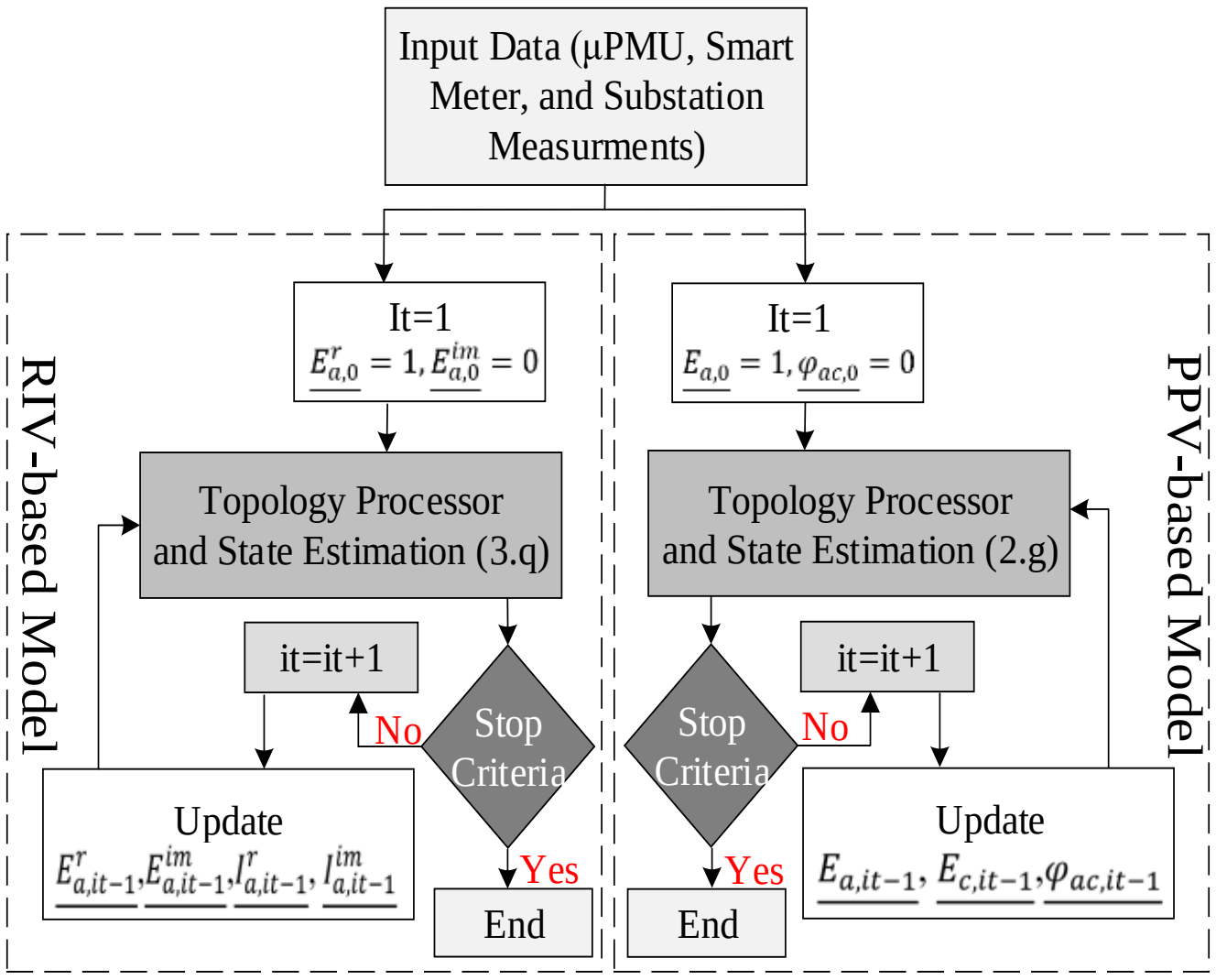

Fig. 1. Proposed iterative MIQP-based topology detection and state estimation in distribution network.

It is worth to note that the accuracy of the proposed linearized models based on first-order approximation of Taylor series is enhanced by solving them iteratively. In the iterative process, the values of $\underline{E_{a,it-1}}$, $\underline{E_{c,it-1}}$, and $\underline{\varphi_{ac,it-1}}$ in the proposed PPV model and $\underline{E^r_{a,it-1}}$, $\underline{E^{im}_{a,it-1}}$, $\underline{I^r_{a,it-1}}$, and $\underline{I^{im}_{a,it-1}}$ in the proposed RIV model are updated using the solution from previous iteration. The simulations for the proposed iterative MIQP-PPV-based and MIQP-RIV-based topology processor and state estimation models are conducted until the stop criteria is met.

## V. Simulation Results

The performances of the proposed MIQP-PPV-based and MIQP-RIV-based topology detection and state estimation methods are demonstrated using a modified IEEE 33-bus distribution system [25]. The test system which is depicted in Fig. 2 includes both radial and meshed topology configurations based on switching actions. The smart meter data are assembled from residential load data of Pecan street Inc. database [26]. For each bus, a random number of houses are selected such that the aggregated load profile of residences follows the nominal value in the IEEE test system. The location and number of micro-PMUs are extracted from [12] and [14] and shown in Fig. 2. In order to calculate actual voltages for various network topologies, nonlinear AC power flow is solved via MATPOWER toolbox in MATLAB [27]. The measurement noise of micro-PMUs is modeled as a Gaussian distribution function with $TVE \leq 0.05\%$ [13]-[14]. The substation injected active and reactive power measurements are also considered, where it is modeled as an ideal voltage source [12]. Smart meters and substation measurements errors are modeled as Gaussian distribution function with errors of 10% and 1%, respectively [10]. To model load's variability, the topology detection and state estimation problem is simulated in 1000 seconds time window with measurement frequency equal to 0.1 $sec^{-1}$, i.e., total of 101 time instants. The proposed model is solved using CPLEX on an Intel Core i7 CPU @ 3.10 GHz computer with 16 GB of RAM.

### A. MIQP-PPV-based Topology Detection and State Estimation

In this section, the performance of the proposed MIQP-PPV-based algorithm in identifying topology of radial and meshed networks is demonstrated by considering the measurement noise of micro-PMUs. Five switches are considered in the test system, which result in $2^5 = 32$ different topologies including radial and meshed configurations. At $x = 440\ sec$, the network topology changes form a radial system with $\Lambda = \{0,0,0,0,0\}$ to a meshed system with $\Lambda = \{1,1,1,0,0\}$ while the status of three switches are changed simultaneously. The simulation is conducted for each time interval and the identified status of switches for the simulated time window is shown in Fig. 3. The results in Fig. 3 confirm that the proposed MIQP-PPV-based topology detection method accurately identifies the radial and meshed topology in all time intervals even while considering load's variability, measurement noise, and multiple simultaneous switching actions. It is worth noting that, the proposed MIQP-PPV-based model detects radial and meshed network topology without knowledge of status of switches, and micro-PMUs and smart meter measurements in prior time intervals.

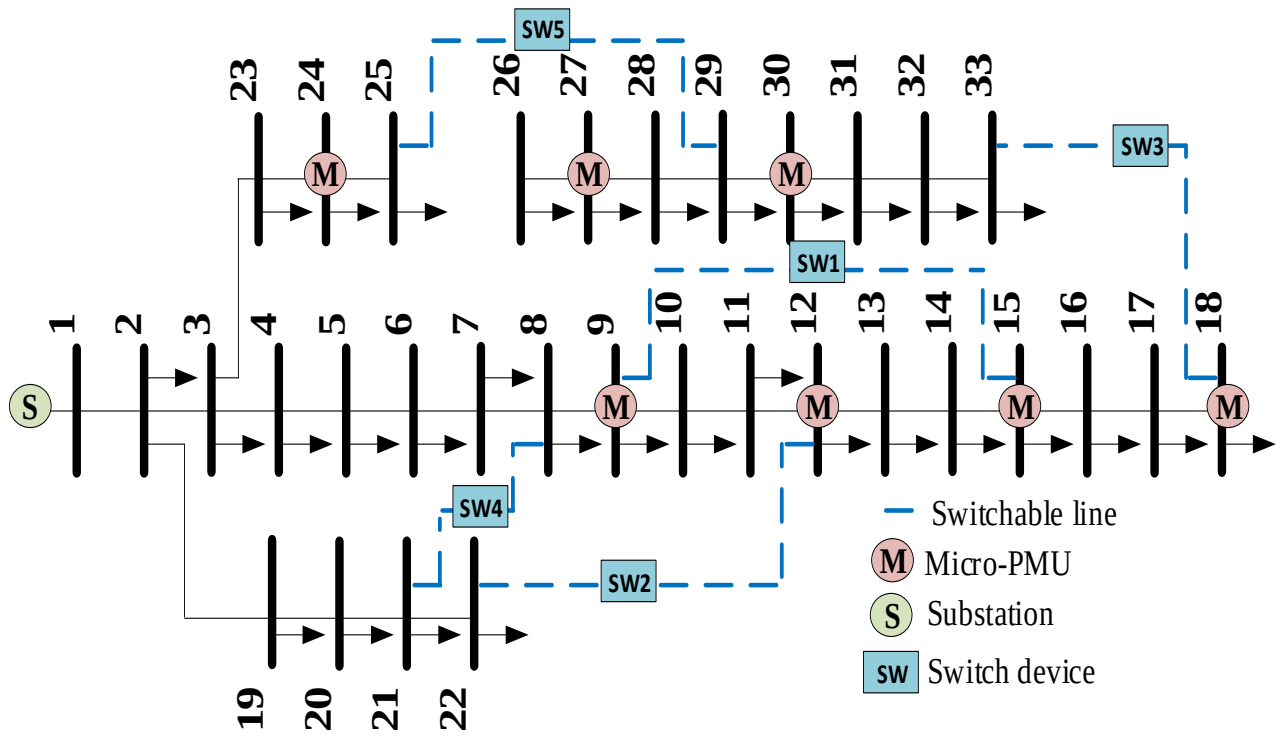

Fig. 2. IEEE 33-bus distribution system equipped with micro-PMUs and switch devices.

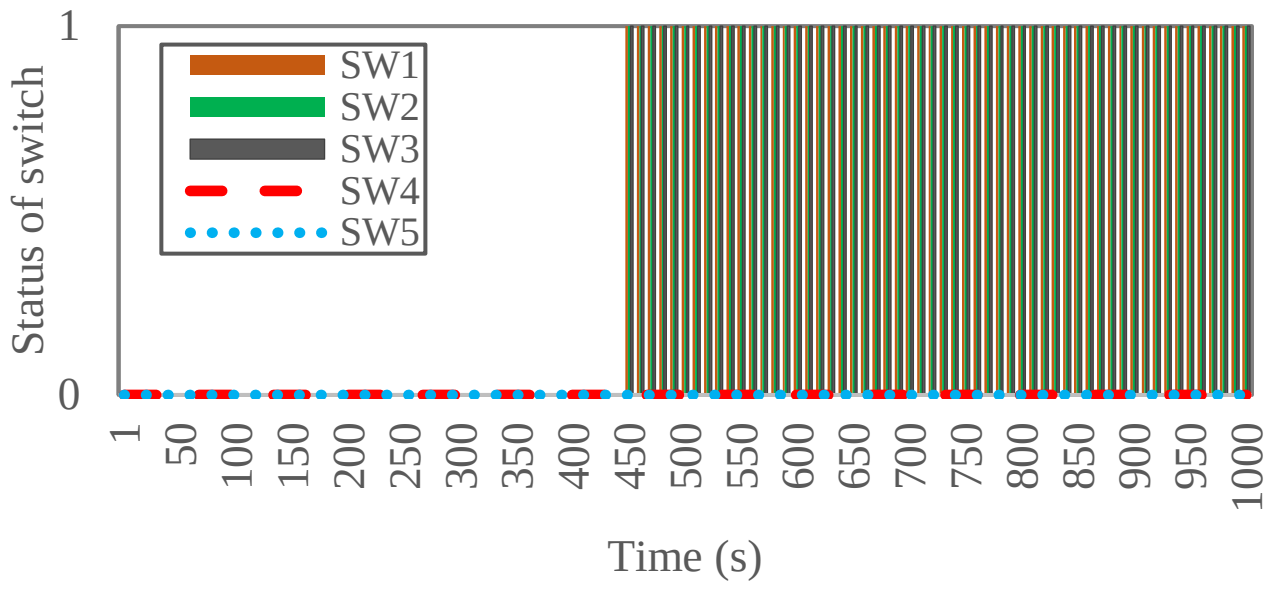

Fig. 3. Status of switches during simulated time window.

Furthermore, the proposed MIQP-PPV-based topology detection model can simultaneously estimate system states. The results of state estimation and the corresponding actual system state values before the topology change, i.e., radial configuration, and after the topology change, i.e., meshed configuration, are compared in Fig. 4 and Fig. 5. These figures confirm that the estimated voltage magnitude and angle closely follow the real voltage profiles in both radial and meshed networks. Also, the absolute error (AE) values of voltage magnitude and error values of voltage angle at each bus are depicted in Figs. 4 and 5. As these figures show, the AE of voltage magnitudes and error of voltage angles for both radial and meshed networks are small. In order to statistically evaluate the performance of the proposed MIQP-PPV-based state estimation model for all time intervals, three indices, namely, root mean square error (RMSE), mean absolute error (MAE), and maximum absolute error (ME) are used. The obtained values of indices for voltage magnitude and angle at each bus over the simulated time window are shown in Figs. 6 and 7. The small values of RMSE, MAE, and ME for all buses confirm that the proposed MIQP-PPV-based model is able to estimate system states with remarkable accuracy.

### B. Method Comparison

In this section, the performance of the proposed MIQP-PPV-based topology processor algorithm is evaluated by comparing with a data-driven method proposed in [12]-[14]. In [14], using the prior information of switch statuses, a library of possible topology configurations based on the change in status of only

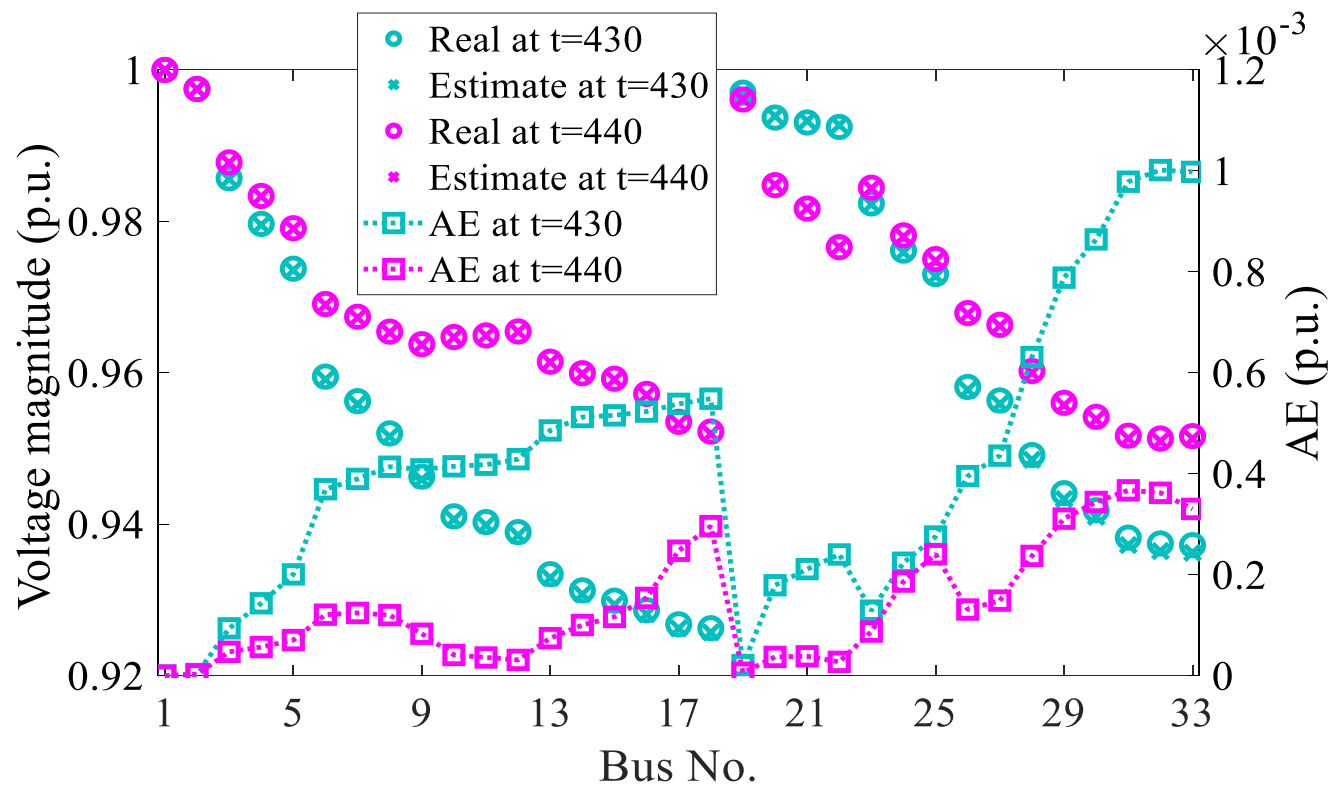

Fig. 4. Estimated, real, and AE of voltage magnitude at x=430 and x=440.

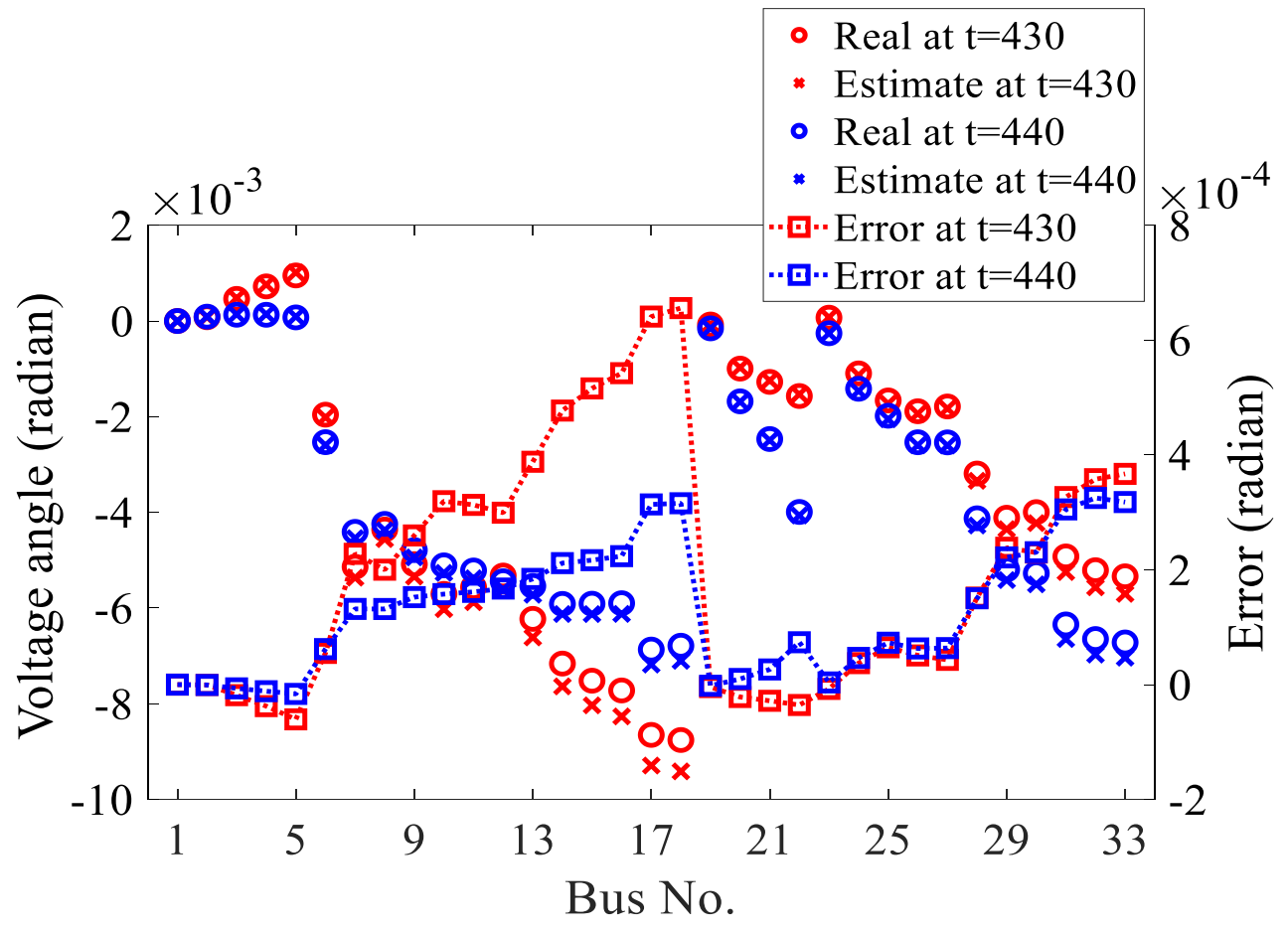

Fig. 5. Estimated, real, and error of voltage angle at x=430 and x=440.

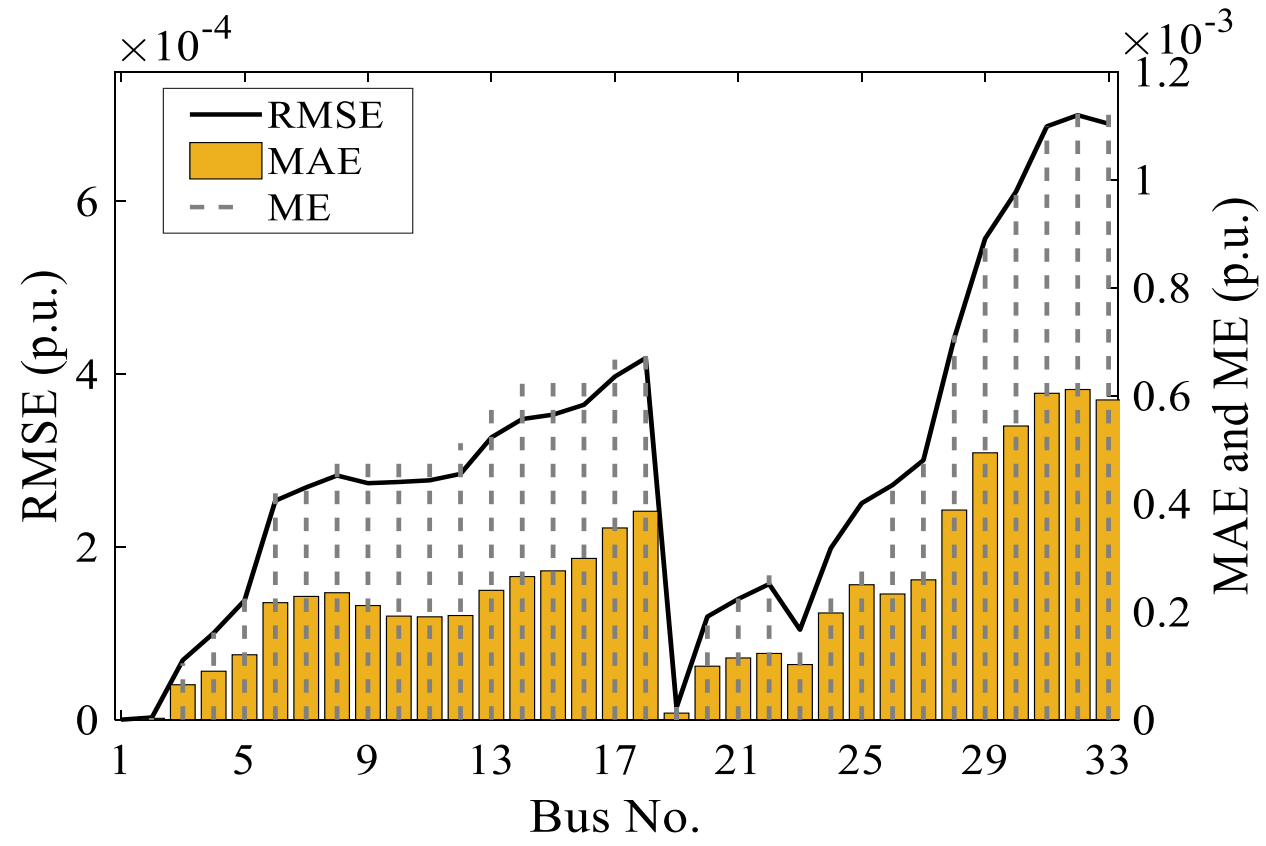

Fig. 6. RMSE, MAE, and ME indices for voltage magnitude at each bus.

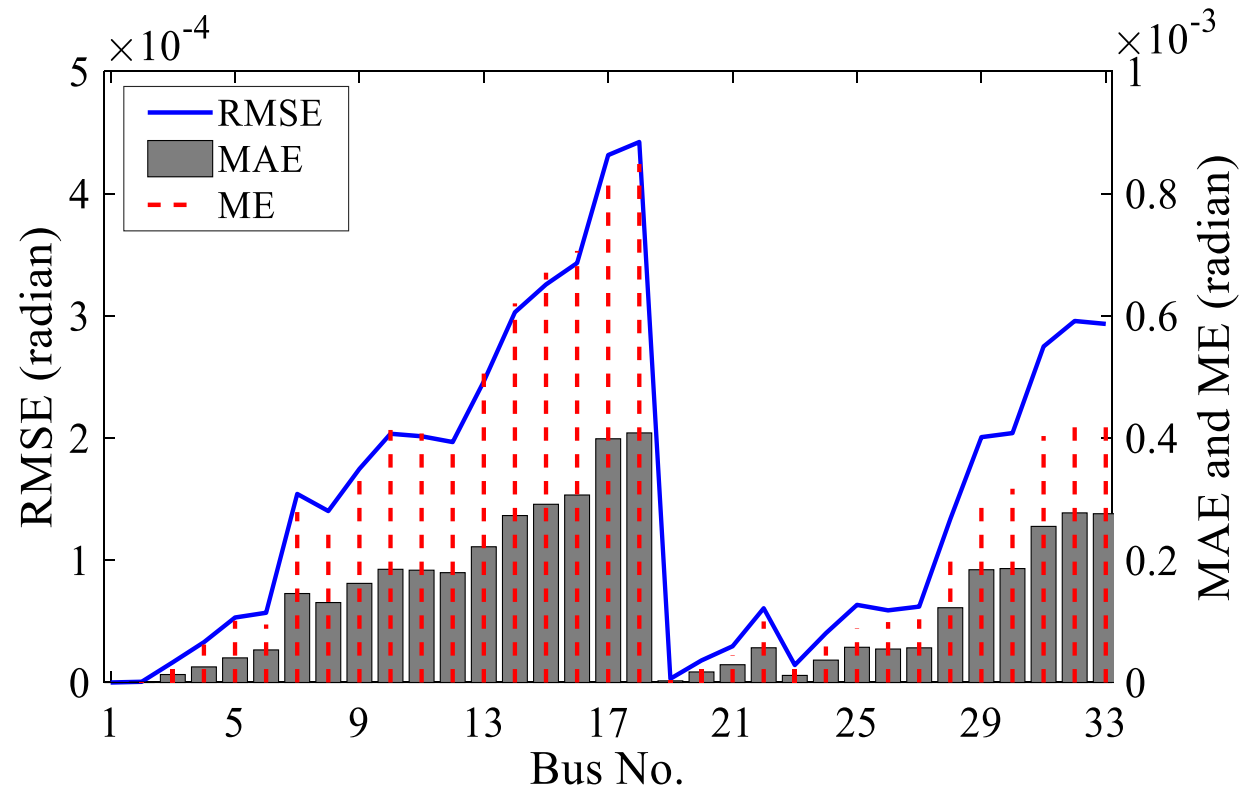

Fig. 7. RMSE, MAE, and ME indices for voltage angle at each bus.

one switch in the system is determined. Then, if the difference between the voltage measured by micro-PMUs at time $x$ (i.e., $E^M_{r,x}$) and time $x-\tau$ (i.e., $E^M_{r,x-\tau}$) is larger than a pre-defined parameter (i.e., min_norm in [14]), it will be projected onto the obtained library of possible system topologies. Finally, the topology with the highest projection value, which is larger than a pre-defined parameter (i.e., min_proj in [14]) is selected as the correct system configuration and topology change time is reported. For the sake of comparison, 100 scenarios are generated based on Monte Carlo simulation while only considering noise for micro-PMUs measurement data. In each scenario, the time interval of topology change within 1000 seconds time window, status of one switch in the system, and measurement noise of micro-PMUs are randomly selected. Four cases are considered for comparing the two methods. In cases 1-3, the smart meter data are collected based on the nominal values of loads provided in IEEE 33-bus test system with different standard deviation (SD) of change of the load between different time intervals. In case 4, the smart meter data are collected from residential load data of Pecan street Inc. database. Tables I-IV compare accuracy of the proposed MIQP-PPV-based topology processor method with the one proposed

TABLE I. COMPARING ACCURACY OF THE PROPOSED MIQP-PPV-BASED METHOD WITH THE METHOD PROPOSED IN [14] WITH SD OF 2.22%.

| **Smart meter data** | **SD** | **Proposed model** | **[14]** | | | |
|---|---|---|---|---|---|---|
| | | ***Accuracy*** | ***min_nor m*** | ***min_proj*** | ***τ*** | ***Accuracy*** |
| IEEE 33-bus | 2.22 | 100% | 0.004 | 0.8 | 5 | 97% |
| | | | 0.006 | 0.9 | 5 | 92% |
| | | | 0.006 | 0.8 | 4 | 91% |

TABLE II. COMPARING ACCURACY OF THE PROPOSED MIQP-PPV-BASED METHOD WITH THE METHOD PROPOSED IN [14] WITH SD OF 3%.

| **Smart meter data** | **SD** | **Proposed model** | **[14]** | | | |
|---|---|---|---|---|---|---|
| | | ***Accuracy*** | ***min_norm*** | ***min_proj*** | ***τ*** | ***Accuracy*** |
| IEEE 33-bus | 3 | 99% | 0.004 | 0.8 | 5 | 82% |
| | | | 0.006 | 0.8 | 5 | 90% |
| | | | 0.008 | 0.9 | 5 | 74% |

TABLE III. COMPARING ACCURACY OF THE PROPOSED MIQP-PPV-BASED METHOD WITH THE METHOD PROPOSED IN [14] WITH SD OF 4%.

| **Smart meter data** | **SD** | **Proposed model** | **[14]** | | | |
|---|---|---|---|---|---|---|
| | | ***Accuracy*** | ***min_norm*** | ***min_proj*** | ***τ*** | ***Accuracy*** |
| IEEE 33-bus | 4 | 97% | 0.004 | 0.8 | 5 | 55% |
| | | | 0.007 | 0.8 | 5 | 79% |
| | | | 0.008 | 0.9 | 5 | 67% |

TABLE IV. COMPARING ACCURACY OF THE PROPOSED MIQP-PPV-BASED METHOD WITH THE METHOD PROPOSED IN [14] USING PECAN STREET DATABASE.

| **Smart meter data** | **Proposed model** | **[14]** | | | |
|---|---|---|---|---|---|
| | ***Accuracy*** | ***min_norm*** | ***min_proj*** | ***τ*** | ***Accuracy*** |
| Pecan Street | 100% | 0.006 | 0.8 | 5 | 49% |
| | | 0.007 | 0.8 | 5 | 58% |
| | | 0.008 | 0.9 | 5 | 57% |

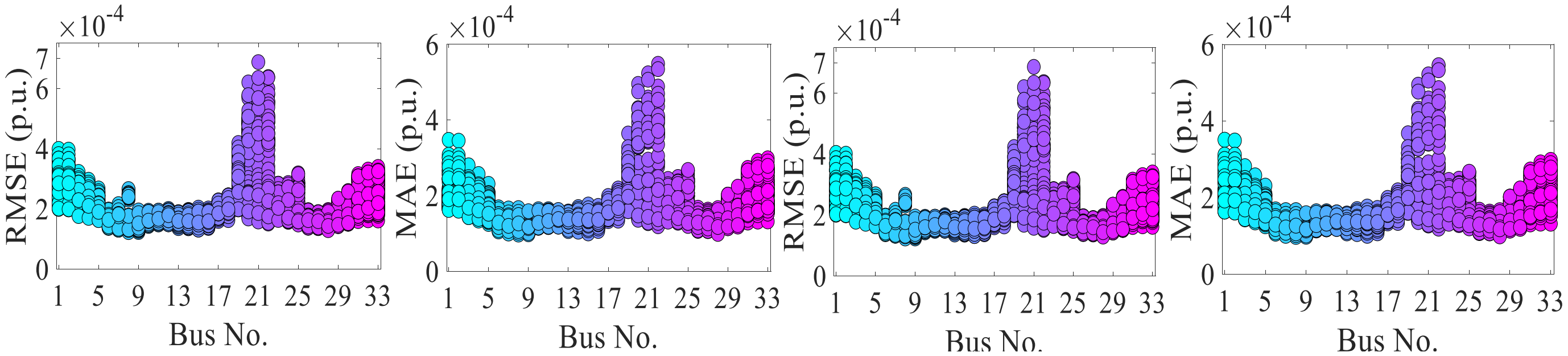

(a) (b) (c) (d)

Fig. 8. Accuracy of voltage magnitude estimation at each bus and scenario for: (a) RMSE and (b) MAE obtained from MIQP-PPV-based model as well as (c) RMSE and (d) MAE obtained from MIQP-RIV-based model.

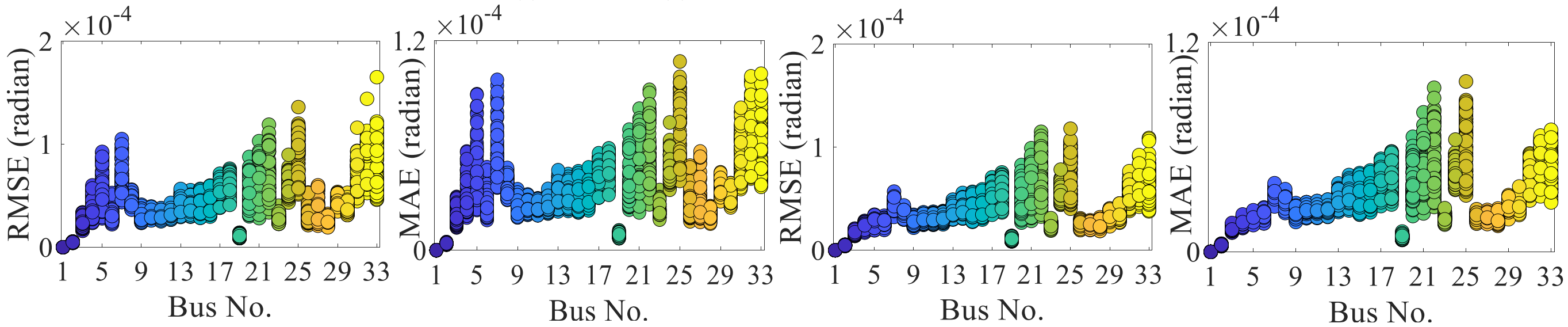

(a) (b) (c) (d)

Fig. 9. Accuracy of voltage angle estimation at each bus and scenario for: (a) RMSE and (b) MAE obtained from MIQP-PPV-based model as well as (c) RMSE and (d) MAE obtained from MIQP-RIV-based model.

in [14] by considering three different parameter tunings for the three parameters (i.e., min_norm, min_proj, and τ) which are used in [14]. According to Table I, the accuracy of the proposed MIQP-PPV-based method among all 100 scenarios is equal to 100% while the accuracy of the model proposed in [14] is dependent to three parameter tunings and at best is equal to 97%. As it can be observed from Tables II and III, by increasing the SD of change of the load, the accuracy of the proposed MIQP-PPV-based method is significantly higher in comparison with the accuracy of the method of [14] with different parameter tunings. In case 4, since the SD of change of the load is high for residential load data of Pecan street Inc. database, the accuracy of the topology detection method proposed in [14] is remarkably low while the proposed MIQP-PPV-based topology processor algorithm identifies the topology with 100% accuracy as shown in Table IV. The reason is that higher load variations, i.e., high SD of change of the load, makes the voltage difference in time series data of micro-PMU measurements to be larger than min_norm parameter; and this change in voltage measurements is projected onto the library of possible system topologies. Therefore, the data-driven method of [14] wrongfully considers the change in the measured voltage time series, which is caused by change of the load, as change in the network topology. Moreover, the method proposed in [14] requires the prior information of switch statuses and measured voltage values by micro-MPUs to identify network topology. In this regard, if prior statuses of switches are wrong, the topology may not be identified correctly. Furthermore, data-driven method in [14] is dependent to three parameter tunings, which limits application of the method in real-time. Since data-driven method in [14] assumes that the topology change may occur due to only one switching action at each time interval, the status of only one random switch in the system is changed at topology transition time in each scenario of Monte Carlo simulation. However, as it is shown in section V-A, the proposed MIQP-PPV-based topology processor model can handle identifying multiple simultaneous switching actions at each time interval without information of switch statuses, micro-PMUs, and smart meters measurements in prior time intervals.

### *C. Comparing performance of proposed MIQP-PPV-based model with proposed MIQP-RIV-based model*

In this section, the performances of the proposed MIQP-PPV-based and the proposed MIQP-RIV-based topology processor and state estimation models are compared by simultaneous modeling of micro-PMUs, smart meters, and substation measurements noise. The simulation is conducted for 100 scenarios, which are generated using Monte Carlo simulation. In each scenario, switches operation time during 1000 seconds time window, status of five switches, and measurement noise of all measurement data are randomly chosen. The accuracy of the proposed MIQP-PPV-based method and the proposed MIQP-RIV-based method among all 100 scenarios with 101 time intervals for the topology identification is 99.83% and 99.84%, respectively.

Since, the proposed models are also able to estimate power system states in the distribution system, the obtained voltage magnitude and angle values from two models are evaluated for each bus and scenario using RMSE and MAE indices as shown in Figs. 8 and 9. The figures confirm that the errors in estimating system states are small with analogous error in terms of voltage magnitude between two methods. However, the MIQP-RIV-based model performs more accurate in terms of estimating voltage angles. The proposed MIQP-RIV-based model outperforms the proposed MIQP-PPV-based model in terms of topology processor and state estimation accuracy. The reason is that in MIQP-RIV-based model, the current flow constraints on the distribution lines are inherently linear, and the only nonlinearity due to inclusion of binary variable associated with status of switches is linearized using big M technique. However, in MIQP-PPV-based model, the AC power flow constraints are linearized in addition to linearization of nonlinearity as a result of adding binary variable associated with status of switches. The average computational time for each snapshot is equal to 0.05 sec and 0.03 sec using the proposed MIQP-PPV-based model and the proposed RIV-based model, respectively, which illustrates the proposed models are computationally efficient for real-time applications.

## VI. Conclusion

In this paper, a simultaneous topology processor and state estimation method is proposed using two mixed-integer quadratic programming (MIQP) formulations, which utilize micro-PMUs and smart meters data. The proposed MIQP approaches are proposed based on two AC optimal power flow models: (i) PPV formulation and (ii) RIV formulation. The results confirm that the proposed MIQP-PPV-based and MIQP-RIV-based models are computationally efficient for real-time application and able to identify different topology configurations including radial and meshed distribution networks. The proposed models are able to detect multiple simultaneous switching actions at each time instant without knowledge of status of switches in prior time intervals. Also, each of the proposed models is a single-shot optimization problem and only requires measurement data at each time snapshot to obtain the topology and states of the system. Monte Carlo simulation is conducted to generate different scenarios of topology and switching actions, switches operation time, and measurement noise. Simulation results illustrate that the proposed models can perform topology identification of a distribution network with high accuracy under load's variability and measurement noises. Moreover, the performance of the proposed MIQP-based state estimation models is examined using statistical indices. The indices confirm that the proposed methods estimate distribution system states with remarkable accuracy. However, the proposed MIQP-RIV-based model outperforms the proposed MIQP-PV-based model in terms of accuracy and speed for topology detection and state estimation.